\documentclass[pdflatex,sn-mathphys-num]{sn-jnl}%

\usepackage{graphicx}%
\usepackage{subcaption}
\usepackage{multirow}%
\usepackage{amsmath,amssymb,amsfonts}%
\usepackage{amsthm}%
\usepackage{mathrsfs}%
\usepackage[title]{appendix}%
\usepackage{xcolor}%
\usepackage{textcomp}%
\usepackage{manyfoot}%
\usepackage{booktabs}%
\usepackage{algorithm}%
\usepackage{algorithmicx}%
\usepackage{algpseudocode}%
\usepackage{listings}%
\usepackage{natbib}
\usepackage{hyperref}

\theoremstyle{thmstyleone}%
\newtheorem{theorem}{Theorem}%
\newtheorem{proposition}[theorem]{Proposition}
\theoremstyle{thmstyletwo}
\newtheorem{example}{Example}%
\newtheorem{remark}{Remark}%
\newtheorem{cor}{Corollary}
\newtheorem{lem}{Lemma}
\theoremstyle{thmstylethree}%
\newtheorem{definition}{Definition}%

\begin{document}

\title[Construction and study of functions with finite energy on the Sierpinski gasket]{Construction and study of functions with finite energy on the Sierpinski gasket}

\author*[1]{\fnm{Harsha} \sur{Gopalakrishnan}}\email{harshagopal3@gmail.com}

\author[2]{\fnm{Srijanani Anurag} \sur{Prasad}}\email{srijanani@iittp.ac.in}

\affil[1]{\orgdiv{Department of Mathematics}, \orgname{SRM Institute of Science and Technology}, \orgaddress{\street{Kattankulathur}, \city{Chennai}, \postcode{603203}, \state{Tamil Nadu}, \country{India}}}

\affil[2]{\orgdiv{Department of Mathematics and Statistics}, \orgname{Indian Institute of Technology  Tirupati}, \orgaddress{\street{Yerpedu–Venkatagiri Road}, \city{Tirupati}, \postcode{517619}, \state{Andhra Pradesh}, \country{India}}}

\abstract{This paper develops a systematic framework for constructing finite-energy functions on the Sierpiński gasket using Lipschitz maps and various classes of contractions. Several properties of the constructed functions are investigated, including the min--max property and bounds for their energy. Several examples of finite-energy functions, along with their graphical representations, are provided to illustrate the functions constructed using the proposed framework. It is further shown that finite energy does not necessarily imply membership in the domain of the Laplacian. The framework is then applied to a class of fractal interpolation functions on the Sierpiński gasket, for which necessary and sufficient conditions for finite energy are established under suitable contractivity assumptions.

}


\maketitle

\section{Introduction}

One of the major areas of research in fractal geometry is analysis on fractals, which originated with Kigami's pioneering work on energy forms and harmonic functions on Sierpiński spaces \cite{kigami1989harmonic}. This theory was subsequently extended to the broader class of post-critically finite (p.c.f.) self-similar sets, which includes the Sierpiński gasket \cite{kigami1993harmonic}. Comprehensive accounts of the subject can be found in the monographs by Kigami \cite{kigami2001analysis} and Strichartz \cite{strichartz2006differential}. The theory of energy forms provides a natural analogue of the classical Dirichlet energy on Euclidean domains and forms the foundation for defining harmonic functions, Laplacians, and other analytical objects on fractals.

Despite the extensive development of analysis on fractals, relatively few explicit classes of functions on the Sierpiński gasket are known to possess finite energy. These include harmonic functions, piecewise harmonic functions, and certain classes of fractal interpolation functions satisfying suitable conditions \cite{ruan2010fractal,sahu2020box,ri2011some,li2017energy}. Furthermore, the domain of the Laplacian is rather restrictive; in particular, it is not closed under pointwise multiplication when the underlying measure is the standard self-similar measure \cite{ben1999not}. Consequently, constructing new classes of finite-energy functions and understanding their analytical properties remain important problems in the analysis on fractals.

Motivated by this, the present paper develops a systematic framework for constructing finite-energy functions on the Sierpiński gasket that extends beyond the classical classes of harmonic and piecewise harmonic functions. Our approach is based on composition with Lipschitz maps and various classes of contractions, providing a unified framework for generating broad families of nonlinear finite-energy functions. We investigate several properties of the constructed functions, including the min--max property, bounds for their energy, and their relationship with the domain of the Laplacian. We provide explicit examples and graphical representations of the constructed functions, illustrating how the proposed framework generates novel classes of functions. We further apply these techniques to fractal interpolation functions and establish necessary and, in the Banach contraction case, sufficient conditions for a  general class of fractal interpolation functions to possess finite energy. For background on the various classes of contractions considered in this paper, we refer the reader to \cite{rakotch1962note,rhoades1977comparison,jachymski2007nonlinear,pasupathi2024very}.

Fractal interpolation on the Sierpiński gasket was initiated by Çelik et al. \cite{ccelik2008fractal} and has subsequently been studied from various perspectives, including interpolation, box dimension, energy, and the min--max property  \cite{ruan2010fractal,sahu2020box,ri2020fractal,de2013pcf, ri2011some}. Since finite energy plays a fundamental role in analysis on fractals, it is natural to investigate conditions under which fractal interpolation functions belong to the energy domain. This question is addressed in this paper. 

The organisation of the paper is as follows: Section~\ref{prelim} presents the necessary preliminaries on analysis on the Sierpiński gasket, including the notions of energy, harmonic functions, various classes of contractions, and Lipschitz maps. Section~\ref{Section1} develops systematic methods for constructing finite-energy functions on the Sierpiński gasket, establishes bounds for their energy, investigates their min--max property, and examines whether these functions belong to the domain of the Laplacian with respect to the standard self-similar measure. Furthermore, various examples and their graphical representations are provided in this section. Finally, Section~\ref{Section2} studies finite-energy fractal interpolation functions on the Sierpiński gasket and establishes necessary and sufficient conditions under which the fractal interpolation functions introduced in \cite{ri2020fractal} possess finite energy.


\section{Some preliminary notions}\label{prelim}

In this section, some preliminary notions and notation that will be used throughout the paper are given. The definitions of various classes of contractions and Lipschitz maps are recalled here. This is followed by a brief review of the basic concepts in the analysis on the Sierpiński gasket, including the standard energy, harmonic and piecewise harmonic functions, and the pointwise Laplacian.

\begin{definition}
    Let $(X,d)$ be a metric space and  $s: X\to X$ be a map.
    \begin{enumerate}
    \item $s$ is called a Banach contraction if there exists a $0\leq c<1$ such that  $d(s(x),s(y))\leq cd(x,y)$ for all $x,y\in X$.
        \item $s$ is called a contractive map if $d(s(x),s(y))<d(x,y)$ for all $x,y\in X$.
        \item $s$ is called a Lipschitz map if there exists a $c>0$ such that  $d(s(x),s(y))\leq cd(x,y)$ for all $x,y\in X$.
    \end{enumerate}
\end{definition}
\begin{definition}\cite{pasupathi2024very}
   Let $(X,d)$ be a metric space, $s:(X,d)\to (X,d)$ and $\phi :\mathbb{R}^+\to \mathbb{R}^+$ be any function. Then,
    \begin{enumerate}
        \item  $s$ is a $\phi$-contraction if
    $d(s(x),s(y))\leq \phi(d(x,y))$ for all $x,y\in X$ and $\phi(t)<t$ for all $t>0$.
       
    \item   $s$ is a Matkowski contraction if  $d(s(x),s(y))\leq \phi(d(x,y))$ for all $x,y\in X$ with nondecreasing $\phi$ and $\lim\limits_{n\to\infty}\phi ^n(t)= 0$ for all $t>0$.

    \end{enumerate}
\end{definition}

\begin{remark}\label{First_rem}
    \begin{enumerate}
        \item Each Banach contraction is a Matkowski contraction, and each Matkowski contraction is a $\phi$-contraction~\cite{pasupathi2024very}.
        \item Each $\phi$-contraction is a contractive map, and each contractive map is a Lipschitz map.
    \end{enumerate}
\end{remark}
The basic notions of analysis on the Sierpiński gasket required in the sequel are recalled below. For a detailed treatment of the standard energy, harmonic and piecewise harmonic functions, and the Laplacian, the reader is referred to \cite{strichartz2006differential}.

Let $V_0=\{P_1,P_2,P_3\}$, where $P_1=(0,0), P_2=(1,0),$ and $P_3=(\frac{1}{2},\frac{\sqrt{3}}{2})$,  denotes the vertices of an equilateral triangle and $u_i:\mathbb{R}^2\to \mathbb{R}^2$ be defined as $u_i(x)=\frac{1}{2}(x+P_i)$ for $x=(x^1,x^2)$ and $i=1,2,3$. The Sierpiński gasket, denoted by $\mathcal{SG}$, is the unique attractor of the iterated function system $\{u_i :i=1,2,3\}$ and satisfies $\mathcal{SG}=u_1(\mathcal{SG})\bigcup u_2(\mathcal{SG})\bigcup u_3(\mathcal{SG})$.

For any word $w=w_1w_2\ldots w_m$ of length $\lvert w\rvert = m$ with $w_i\in\{1,2,3\}$ for $i=1,2,\ldots ,m$, define $u_w=u_{w_1}\circ u_{w_2}\circ \ldots\circ u_{w_m}$ and $V_m=\bigcup\limits_{\lvert w \rvert=m}u_w(V_0)$. Then, $\mathcal{SG}$ can also be defined as $\mathcal{SG}=\overline{\bigcup\limits_{m=0}^{\infty}V_m}$. Let $\Gamma_0$ be the complete graph on $V_0$. The graph $\Gamma_m$ on $V_m$ is defined inductively as follows: Two points $x,y\in V_m$ are adjacent to each other, denote $x\sim_m y$, if and only if $x=u_i(x')$, $y=u_i(y')$ for some $u_i$, and $x'$ and $y'$ are adjacent in $\Gamma_{m-1}$. 

The standard energy is the fundamental quadratic form in the analysis on the Sierpiński gasket. It plays a role analogous to the classical Dirichlet energy $\int |\nabla u|^2$ on Euclidean domains and forms the basis for defining harmonic functions, the Laplacian, and the associated energy space on $\mathcal{SG}$. Although this paper focuses exclusively on the standard energy, a bilateral energy form corresponding to each $b>0$, can also be defined on the Sierpiński gasket. Interested readers are referred to \cite{strichartz2006differential,kigami2001analysis} for a brief introduction and to \cite{sabot1997existence} for a more detailed study in this direction. We recall the standard energy and the related concepts that will be used throughout the paper. The definitions required in the sequel are given below.

\begin{definition}\cite{kigami2001analysis,strichartz2006differential}
For any $f:\mathcal{SG}\to \mathbb{R}$, the standard energy of $f$ is given by $\varepsilon(f)=\lim\limits_{m\to\infty}\varepsilon_m(f)$, where $\varepsilon_m(f)=\big(\frac{5}{3}\big)^m\sum\limits_{x\sim_m y}^{}(f(x)-f(y))^2$ is the $m^{th}$ stage energy of $f$.
\end{definition}
\begin{definition}\cite{kigami2001analysis,strichartz2006differential} \label{definition}
For $f,g:\mathcal{SG}\to \mathbb{R}$, the standard bilinear energy form of the functions $f$ and $g$ is given by $\varepsilon(f,g)=\lim\limits_{m\to\infty}\varepsilon_m(f,g)$, where $\varepsilon_m(f,g)=\big(\frac{5}{3}\big)^m\sum\limits_{x\sim_m y}^{}(f(x)-f(y))(g(x)-g(y))$. Note that $\varepsilon(f,f)=\varepsilon(f)$.
\end{definition}
\begin{definition}\cite{kigami2001analysis,strichartz2006differential}
A function $h:\mathcal{SG}\to\mathbb{R}$ is a harmonic function if $\varepsilon_m(h)=\varepsilon_{m-1}(h)$ for any $m\geq 1$.
\end{definition}
\begin{definition}\cite{kigami2001analysis,strichartz2006differential}
A function $h':\mathcal{SG}\to \mathbb{R}$ is a piecewise harmonic function if there exist an $m'\geq 1$ such that $\varepsilon_m(h')=\varepsilon_{m-1}(h')$ for all $m\geq m'$. The case $m'=1$ corresponds to a harmonic function.
\end{definition}
Let the function values at the boundary set $V_0$ be given by $h(P_1),h(P_2)$ and $h(P_3)$.  Given these boundary values, a unique harmonic extension to $V_1$ is obtained by means of the $\frac25-\frac15$ rule: $h(P_{12})=\frac{2}{5}h(P_1)+\frac{2}{5}h(P_2)+\frac{1}{5}h(P_3)$, $h(P_{13})=\frac{2}{5}h(P_1)+\frac{1}{5}h(P_2)+\frac{2}{5}h(P_3)$ and $h(P_{23})=\frac{1}{5}h(P_1)+\frac{2}{5}h(P_2)+\frac{2}{5}h(P_3)$, where $P_{ij}$ denotes the midpoint of the line segment joining $P_i$ and $P_j$ for $i,j\in\{1,2,3\}$ and $i\neq j$. Extend $h$ to $V_m$ using the $\frac{2}{5}-\frac{1}{5}$ rule for any $m\geq 1$ iteratively and further extend it continuously to $\mathcal{SG}$, in order to obtain the harmonic function on $\mathcal{SG}$. The collection of all harmonic functions on $\mathcal{SG}$ forms a vector space of dimension $3$.

Let $C(\mathcal{SG})=\{f:\mathcal{SG}\to\mathbb{R}:f\text{ is continuous}\}$ and $\mathrm{dom}(\varepsilon)=\{f\in C(\mathcal{SG}):\varepsilon(f)~<~\infty\}.$  The space $\mathrm{dom}(\varepsilon)$ is called the energy domain of the Sierpiński gasket, and its elements are referred to as finite-energy functions. It is well known that $\mathrm{dom}(\varepsilon)$ is dense in $C(\mathcal{SG})$. Moreover, $\mathrm{dom}(\varepsilon)$ modulo constants forms a real inner product space with inner product $\varepsilon(f,g)$. It is also known that $\mathrm{dom}(\varepsilon)$ is closed under pointwise multiplication.

Next, we will recall the pointwise definition of the Laplacian together with a well-known result concerning its domain.

\begin{definition}\cite{kigami2001analysis,strichartz2006differential}
For $f\in\mathrm{dom}(\varepsilon)$, and $x\in \mathcal{SG}$, the pointwise Laplacian at $x$ is defined to be 
\begin{equation}\label{lapla}
\Delta f(x)=\frac{3}{2}\lim\limits_{m\to\infty} 5^m\Delta_m f(x),
\end{equation}
where $\Delta_m f(x)=\sum\limits_{y\sim_m x}(f(y)-f(x))$.
The function $f$ is said to belong to the domain of the Laplacian, denoted by $\mathrm{dom}(\Delta)$, if and only if the limit in Eq. \eqref{lapla} exists uniformly.
\end{definition}

\begin{proposition}\cite{strichartz2006differential}\label{whatnot}
If $f$ and $g$ are two non-constant functions in $\mathrm{dom}(\Delta)$, then $fg$ does not belong to $\mathrm{dom}(\Delta)$.
\end{proposition}
The above proposition shows that the domain of the Laplacian is considerably more restrictive than the energy domain. This distinction will be used in Section~\ref{Section1}. 
\vspace{0.2cm}
\par
This section is concluded by briefly introducing the concept of fractal interpolation on the Sierpinski gasket. The pioneering work in this direction was initiated by Çelik et al.~\cite{ccelik2008fractal}. The idea of fractal interpolation on the Sierpinski gasket is to construct an analogue of Barnsley’s fractal interpolation on intervals \cite{barnsley1988fractals}.
The main result on fractal interpolation on the Sierpinski gasket from \cite{ccelik2008fractal} is provided as follows.
\begin{theorem}\cite{ccelik2008fractal}\label{Intrpoln_thm1}
    Let $\nu: V_n\to \mathbb{R}$ be any function, which can be considered as the set of interpolation values on the set $V_n$, and $\alpha_w$ be any given numbers with $0<|\alpha_w|<1$ for each word $w$ of length $n$. Then, there exists a unique continuous function $f:\mathcal{SG}\to\mathbb{R}$, such that $f|_{V_n}=\nu$ and
    \begin{equation*}
        f(u_w(x))=\alpha_w f(x)+h_w(x) \text{ for } x\in\mathcal{SG},
    \end{equation*}
    where $h_w$ are harmonic functions with boundary values
    \begin{equation*}
        h_w(P_i)=\nu(u_w(P_i))-\alpha_w\nu(P_i) \text{ for each } i=1,2,3.
    \end{equation*}
\end{theorem}
Here, for any initial function $\nu$ on $V_n$, a self-similar interpolation on $\mathcal{SG}$ is obtained with the use of harmonic functions. A general version of this result with $n=1$ is given in \cite{ri2020fractal}, and a brief description of the same will be explained in Section \ref{Section2}.

\section{New methods to construct finite energy functions}\label{Section1}
This section begins with the definition of the min--max property, followed by the construction of new classes of finite-energy functions and the study of their properties. 

\begin{definition}
Let $f:\mathcal{SG}\to\mathbb{R}$ be a function in $\mathrm{dom}(\varepsilon)$. The function $f$ is said to satisfy the \emph{min--max property} if
\[
\min_{1\le i\le3}f(P_i)\le f(x)\le \max_{1\le i\le3}f(P_i), \qquad x\in\mathcal{SG}.
\]
Equivalently, the minimum and maximum values of $f$ are attained on the boundary set $V_0$.
\end{definition}

\begin{remark}
It is well known that every harmonic function on $\mathcal{SG}$ satisfies the min-max property.
\end{remark}
The following proposition is immediate from the definition of the min--max property, and its proof is therefore omitted.
\begin{proposition}
If $h$ is a harmonic function, then $\frac{1}{h}$ satisfies the min-max property, provided $h(x)\neq 0$ for all $x\in \mathcal{SG}$. Moreover, if $f\in \mathrm{dom}(\varepsilon)$ satisfies the min-max property with $f(x)\neq 0$ for all $x\in \mathcal{SG}$, then $\frac{1}{f}$ also satisfies the min-max property.
\end{proposition}

The next result shows that the reciprocal of a non-vanishing harmonic function has finite energy.
 
\begin{theorem}\label{first_thm}
If $h$ is a harmonic function with $h(x)\neq 0$ for all $x\in\mathcal{SG}$, then the function $\frac{1}{h}:\mathcal{SG}\to \mathbb{R}$ defined by $\frac{1}{h}(x)=\frac{1}{h(x)}$ is in $\mathrm{dom}(\varepsilon)$.
\end{theorem}
\begin{proof}
The energy of the function $\frac{1}{h}$ is,
\begin{align*}
\varepsilon\left(\frac{1}{h}\right)=\lim\limits_{m\to\infty}\left(\frac{5}{3}\right)^m\sum\limits_{x\sim_m y}^{}\left(\frac{1}{h(x)}-\frac{1}{h(y)}\right)^2=\lim\limits_{m\to\infty}\left(\frac{5}{3}\right)^m\sum\limits_{x\sim_m y}^{}\frac{\left(h(y)-h(x)\right)^2}{h(x)^2h(y)^2} \ \leq c\varepsilon(h) 
\end{align*}
 where, $c=\left(\min\limits_{1\leq i\leq 3} h(P_i)\right)^{-4}$.  Since harmonic functions satisfy the min-max property, $h$ is bounded away from zero on $\mathcal{SG}$, and therefore $c<\infty$. Hence, $\varepsilon\!\left(\frac1h\right)<\infty,$ which proves the theorem.  
 \end{proof}
 \begin{remark}\label{Rem_for1over_h_first}
    Note that the function $\frac{1}{h}$ given in Theorem \ref{first_thm} is not a harmonic function when $h$ is a non-constant harmonic function. This is because it does not satisfy the harmonic extension algorithm. The detailed description of this fact can be seen as a consequence of the proof of claim 3 in Theorem \ref{harm_decomp}. More details are provided in Remark \ref{Rem_for1over_h}. Thus, a non-harmonic function is obtained which nevertheless has finite energy and admits an explicit description. A detailed description of this kind of function is provided through the following example.
 \end{remark}
\begin{example}\label{EX1}
   Let $h_1$ be the harmonic function with boundary values $h_1(P_1)=2$, $h_1(P_2)=5$, and $h_1(P_3)=1$. Hence, by the harmonic extension algorithm, it follows that $h_1(P_{12})=3$, $h_1(P_{13})=\frac{11}{5}$, and $h_1(P_{23})=\frac{14}{5}$. Thus,
    $\frac{1}{h_1}(P_1)=\frac{1}{2},
    \frac{1}{h_1}(P_2)=\frac{1}{5},\frac{1}{h_1}(P_3)=1,\frac{1}{h_1}(P_{12})=\frac{1}{3},\frac{1}{h_1}(P_{13})=\frac{5}{11},$ and $\frac{1}{h_1}(P_{23})=\frac{5}{14}$. Clearly, $\frac{1}{h_1}$ is not a harmonic function, but it is a finite energy function, since $\varepsilon(\frac{1}{h_1})\leq \varepsilon(h_1)=26$ by Theorem \ref{first_thm}. A pictorial representation of the functions $h_1$ and $\frac{1}{h_1}$ is given in Fig. \ref{fig:h_function} and Fig. \ref{fig:inverse_h_function} respectively. 
    
    Note that the scaling of the $z$-axis differs between the two figures. This difference in scaling is intentional and is used to provide a better visual representation of the two functions.

\begin{figure}[htbp]
    \centering
    \begin{minipage}{0.48\textwidth}
        \centering
        \includegraphics[width=\textwidth]{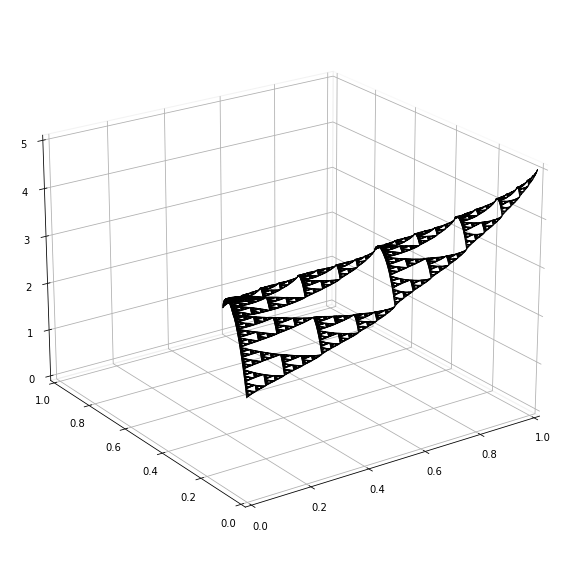}
        \caption{Pictorial representation of the harmonic function $h$ with boundary values $h_1(P_1)=2$, $h_1(P_2)=5$, and $h_1(P_3)=1$.}
        \label{fig:h_function}
    \end{minipage}
    \hfill
    \begin{minipage}{0.48\textwidth}
        \centering
        \includegraphics[width=\textwidth]{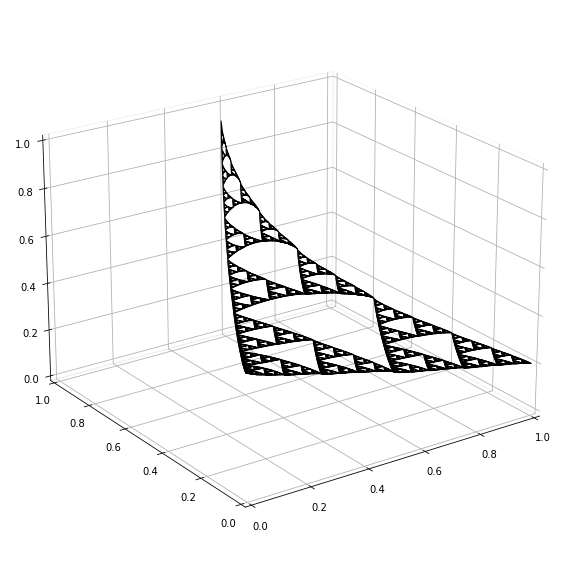}
        \caption{Pictorial representation of the finite energy function $\frac{1}{h_1}$.}
        \label{fig:inverse_h_function}
    \end{minipage}
\end{figure}
    
\end{example}

 The following theorem provides a general method for constructing finite-energy functions from existing ones.
\begin{theorem}\label{first lemma}
If $f\in \mathrm{dom}(\varepsilon)$ and $s:\mathbb R\to\mathbb R$ is a Lipschitz map, then $s\circ f\in \mathrm{dom}(\varepsilon)$.
\end{theorem}
\begin{proof}
The $m^{th}$ stage energy of $s\circ f$ is
\begin{align*}
\varepsilon_m(s\circ f)=\Big(\frac{5}{3}\Big)^m\sum\limits_{x\sim_m y}
\lvert s\circ f(x)-s\circ f(y)\rvert^2&=\Big(\frac{5}{3}\Big)^m\sum\limits_{x\sim_m y}
d( s\circ f(x), s\circ f(y))^2\\&  =\Big(\frac{5}{3}\Big)^m\sum\limits_{x\sim_m y}  d( s(f(x)), s(f(y)))^2.
\end{align*} 
Since $s$ is a Lipschitz map, there exists a $c>0$ such that 
\begin{align*}
d(s(f(x)), s(f(y)))\leq c d(f(x),f(y)).
\end{align*}
Hence, 
\begin{align*}
\varepsilon_m(s\circ f)\leq c^2 \Big(\frac{5}{3}\Big)^m\sum\limits_{x\sim_m y} d(f(x),f(y))^2 &=c^2\Big(\frac{5}{3}\Big)^m\sum\limits_{x\sim_m y}\lvert f(x)-f(y)\rvert^2\\&=c^2\varepsilon_m(f) \quad \leq c^2\varepsilon(f) \quad <\infty.
\end{align*}
Thus,
\[
\varepsilon(s\circ f) = \lim_{m\to\infty}\varepsilon_m(s\circ f) \leq c^2\varepsilon(f) < \infty,
\]
and hence $s\circ f\in \mathrm{dom}(\varepsilon)$.
\end{proof}
\begin{example}\label{EX2}
Let $s:[0,\frac{\pi}{2}]\to[0,1]$ be defined as $s(x)=\sin(x)$. Then, $s$ is continuous on $[0,\frac{\pi}{2}]$ and differentiable on $(0,\frac{\pi}{2})$. By the mean value theorem, it is easy to see that $|s(x)-s(y)|\leq |x-y|$ for all $x,y\in[0,\frac{\pi}{2}]$ and so $s$ is a Lipschitz map.
Consider the harmonic function $h_2$ with boundary values $h_2(P_1)=0, h_2(P_2)=\frac{\pi}{2},$ and $h_2(P_3)=\frac{\pi}{2}$. Then, $s\circ h_2:\mathcal{SG}\to [0,1]$ given by $s\circ h_2(x)=\sin (h_2(x))$ is not a harmonic function, but a finite energy function. The boundary values of $s\circ h_2$ are given by $s\circ h_2(P_1)=0, s\circ h_2(P_2)=1,$ and $s\circ h_2(P_3)=1$. Furthermore, $2<\varepsilon(s\circ h_2)\leq \varepsilon(h_2)=\frac{\pi^2}{4}\approx 2.4674$, where the lower bound is the energy of the harmonic function with boundary values $0,1$, and $1$ at $P_1, P_2$, and $P_3$ respectively. The pictorial representation of $h_2$ and $s\circ h_2$ is given in Fig. \ref{fig:h2_function} and Fig. \ref{fig:sinh2_function}, respectively.

\begin{figure}[htbp]
\centering
\begin{minipage}{0.48\textwidth}
\centering
\includegraphics[width=\textwidth]{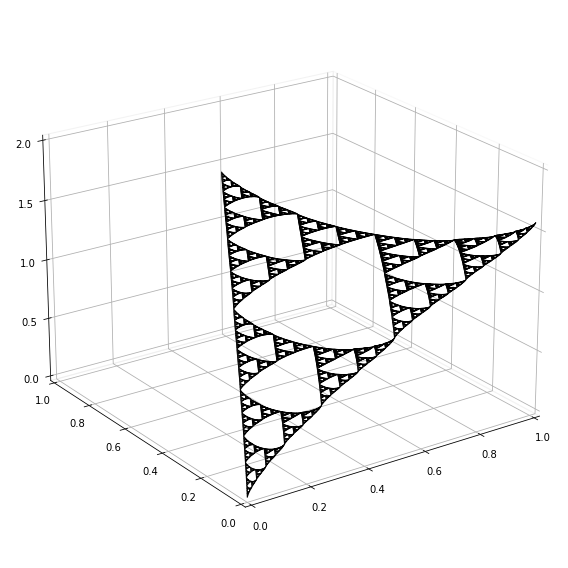}
\caption{Pictorial representation of the harmonic function $h_1$ with boundary values $h_2(P_1)=0$, $h_2(P_2)=\frac{\pi}{2}$, and $h_3(P_3)=\frac{\pi}{2}$.}
        \label{fig:h2_function}
    \end{minipage}
    \hfill
    \begin{minipage}{0.48\textwidth}
        \centering
        \includegraphics[width=\textwidth]{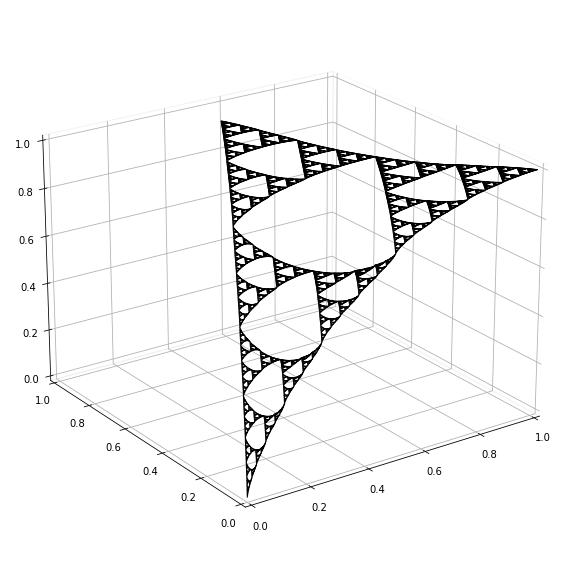}
        \caption{Pictorial representation of the finite energy function $s\circ h_2$.}
        \label{fig:sinh2_function}
    \end{minipage}
\end{figure}
\end{example}
\begin{remark}

The converse of Theorem~\ref{first lemma} does not hold in general, even when $s$ is a Banach contraction. The following examples illustrate both possibilities.

\begin{enumerate}
\item Let $s$ be a constant, then $s\circ f$ is a constant and hence $s\circ f$ is a harmonic function. Thus, $s\circ f\in \mathrm{dom}(\varepsilon)$.
\item Let $s=\frac{x}{2}$. Then, 
\begin{equation}\nonumber
\varepsilon_m(s\circ f) =\Big(\frac{5}{3}\Big)^m\frac{1}{4}\sum\limits_{x\sim_m y} \lvert f(x)-f(y)\rvert^2=\frac{1}{4}\varepsilon_m(f).
\end{equation}
Since $f$ is of infinite energy, $\varepsilon_m(f)\to\infty$ as $m\to\infty$. Thus, $\varepsilon_m(s\circ f)\to\infty$ as $m\to\infty$ giving $s\circ f$ is of infinite energy.
\end{enumerate}
\end{remark}
\begin{cor}\label{main_theorem}
Let $f$ be any function in $\mathrm{dom}(\varepsilon)$ with $\displaystyle{f(x)\in\mathbb{R}^+\cup\{0\}}$ for all $x\in\mathcal{SG}$. Consider the functions $g_1,g_2:\mathcal{SG}\to \mathbb{R}$ as $\displaystyle{g_1(x)=\frac{1}{1+f(x)}}$ and $\displaystyle{g_2(x)=\frac{f(x)}{1+f(x)}}$. Then, $g_1,g_2\in \mathrm{dom}(\varepsilon)$. 
\end{cor}
\begin{proof}
Let $s_1,s_2:\mathbb{R}^+\cup\{0\}\to \mathbb{R}^+\cup\{0\}$ be defined by $s_1(t)=\frac{1}{1+t}$ and $s_2(t)=\frac{t}{1+t}$. Both $s_1$ and $s_2$ are Lipschitz maps. Therefore, Theorem~\ref{first lemma} immediately yields $g_1=s_1\circ f,$ and $ g_2=s_2\circ f,$ are functions on $\mathcal{SG}$ with finite energy and hence $g_1,g_2\in\mathrm{dom}(\varepsilon).$ 
\end{proof}  

In particular, the preceding corollary also holds for harmonic functions. However, a stronger result for harmonic functions follows from Theorem \ref{first_thm}.

\begin{cor}\label{coro_har}
If $h$ is a harmonic function with $h(x)\neq -1$ for all $x\in\mathcal{SG}$, then $\frac{1}{1+h}$ and $\frac{h}{1+h}$ are in $\mathrm{dom}(\varepsilon)$.
\end{cor}
\begin{proof}
Since $1+h$ is a harmonic function, Theorem \ref{first_thm} gives  $\frac{1}{1+h}\in \mathrm{dom}(\varepsilon)$. Also, since $\frac{h}{1+h}=1-\frac{1}{1+h}$ and both the constant function $1$ and $\frac{1}{1+h}$ have finite energy, $\frac{h}{1+h}\in \mathrm{dom}(\varepsilon)$.
\end{proof}

The properties of the functions obtained in Corollary \ref{main_theorem} are now investigated.

\begin{theorem}
Let $f$ be any function in $\mathrm{dom}(\varepsilon)$ with $\displaystyle{f(x)\in\mathbb{R}^+\cup\{0\}}$ for all $x\in\mathcal{SG}$ and satisfies the min-max property. Then, $\frac{1}{1+f}$ and $\frac{f}{1+f}$ satisfy the min-max property.
\end{theorem}
\begin{proof}
Since $f$ satisfies the min--max property,
\[
\frac{1}{1+\max_{1\le i\le3}f(P_i)} \leq \frac{1}{1+f(x)} \leq \frac{1}{1+\min_{1\le i\le3}f(P_i)},
\]
for every $x\in\mathcal{SG}$. Hence, $\frac{1}{1+f}$ satisfies the min--max property. Furthermore, $\frac{f}{1+f} = 1-\frac{1}{1+f}.$ Since both the constant function $1$ and the function $\frac{1}{1+f}$ satisfy the min-max property, the function $\frac{f}{1+f}$ also satisfies the min--max property.
\end{proof}

The following theorem establishes some basic energy properties of the functions $\frac{1}{1+f}$ and $\frac{f}{1+f}$ for $f\in \mathrm{dom}(\varepsilon)$.

\begin{theorem}
Let $f\in \mathrm{dom}(\varepsilon)$ such that $\displaystyle{f(x)\in\mathbb{R}^+\cup\{0\}}$ for all $x\in\mathcal{SG}$. Then,
\begin{enumerate}
\item $\varepsilon\left(\frac{1}{1+f}\right)=\varepsilon\left(\frac{f}{1+f}\right)$.
\item  $\varepsilon(h_1^*)=\varepsilon(h_2^*)$, where $h_1^*$ and $h_2^*$ are the unique harmonic functions with boundary values $\frac{1}{1+f(P_i)}$ and $\frac{f(P_i)}{1+f(P_i)}$ respectively for $i=1,2,3$. 
\item  $\varepsilon(h^*)\leq \varepsilon\left(\frac{1}{1+f}\right)=\varepsilon\left(\frac{f}{1+f}\right)\leq \varepsilon(f)$, where $\varepsilon(h^*)=\varepsilon(h_1^*)=\varepsilon(h_2^*)$.
\end{enumerate}
\end{theorem}
\begin{proof}
\begin{enumerate}
\item Note that $\varepsilon(1,v)=0$ for all $v\in \mathrm{dom}(\varepsilon)$. Hence $\varepsilon\left(\frac{1}{1+f}+\frac{f}{1+f},v\right)=\varepsilon(1,v)=0$ for all $v\in \mathrm{dom}(\varepsilon)$. In particular, when $v=\frac{f}{1+f}$,
\begin{equation}\label{eq_1}
\varepsilon\left(\frac{f}{1+f}\right)=-\varepsilon\left(\frac{1}{1+f},\frac{f}{1+f}\right)
\end{equation}
and when $v=\frac{1}{1+f}$,
\begin{equation}\label{eq_2}
\varepsilon\left(\frac{1}{1+f}\right)=-\varepsilon\left(\frac{f}{1+f},\frac{1}{1+f}\right).
\end{equation}
Combining the Eqs. \eqref{eq_1} and \eqref{eq_2}, the required result is obtained.

\item Since $h_1^*$ and $h_2^*$ are harmonic functions, $\varepsilon(h_1^*)$ and $\varepsilon(h_2^*)$ coincide with their respective zeroth-stage energies. Hence, it is easy to see that $\varepsilon(h_1^*)=\varepsilon(h_2^*)$.
\item Since a harmonic function minimizes the energy among all functions with the same boundary values, 
\begin{align*}
    \varepsilon(h^*)\leq \varepsilon\left(\frac{1}{1+f}\right)=\varepsilon\left(\frac{f}{1+f}\right).
\end{align*} Also, by Theorem \ref{first lemma}
and Corollary \ref{main_theorem}, it follows that $\varepsilon\left(\frac{1}{1+f}\right)=\varepsilon\left(\frac{f}{1+f}\right)\leq \varepsilon(f)$.
\end{enumerate}
\end{proof}
\begin{example}\label{EX3}
Consider the harmonic function with boundary values $h_3(P_1)=0, h_3(P_2)=0,$ and $h_3(P_3)=1$. The pictorial representations of $h_3,\frac{1}{1+h_3},$ and $\frac{h_3}{1+h_3}$ are given in Fig.~\ref{fig:three_figures}. 
Furthermore, $\frac{1}{2}\leq \varepsilon(\frac{1}{1+h_3})=\varepsilon(\frac{h_3}{1+h_3})\leq 2$.
 \begin{figure}[htbp]
    \centering
    \begin{subfigure}{0.65\textwidth}
        \centering
        \includegraphics[width=\linewidth]{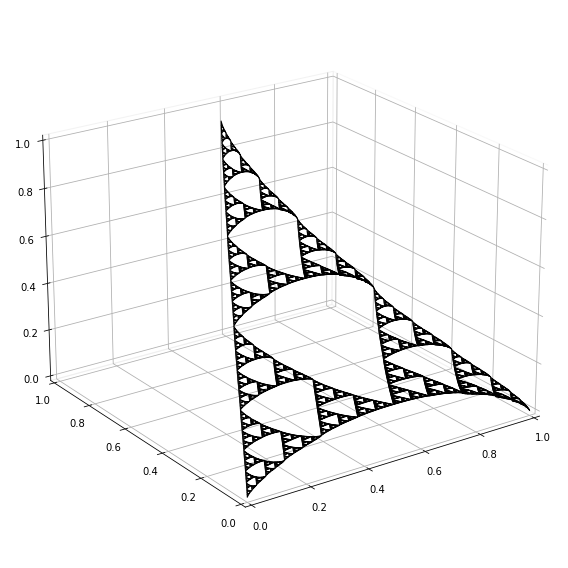}
        \caption{Harmonic function $h_3$ with boundary values $h_3(P_1)=0, h_3(P_2)=0,$ and $h_3(P_3)=1$.}
        \label{fig:1a}
    \end{subfigure}
    \vspace{0.6cm}
    \begin{subfigure}{0.48\textwidth}
        \centering
        \includegraphics[width=\linewidth]{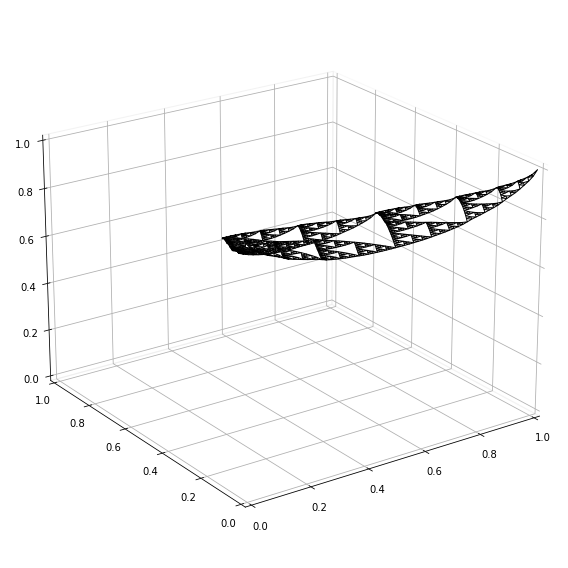}
        \caption{The function $\frac{1}{1+h_3}$}
        \label{fig:1b}
    \end{subfigure}
    \hfill
    \begin{subfigure}{0.48\textwidth}
        \centering
        \includegraphics[width=\linewidth]{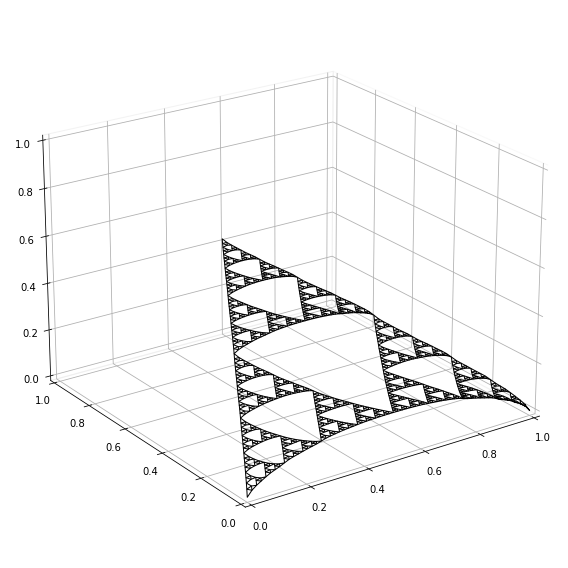}
        \caption{The function $\frac{h_3}{1+h_3}$}
        \label{fig:1c}
    \end{subfigure}
    \caption{Pictorial representation of $h_3, \frac{1}{1+h_3}$ and $\frac{h_3}{1+h_3}$}
\label{fig:three_figures}
\end{figure}
\end{example}
\begin{remark}
 Examples \ref{EX1},\ref{EX2}, and \ref{EX3}  illustrate that our new results yield a collection of finite-energy functions. 
 In the same way, in Example \ref{EX3}, instead of a harmonic function, we can choose any function
$f\in\operatorname{dom}(\varepsilon)$ to construct new functions in
$\operatorname{dom}(\varepsilon)$. For example, consider
$f=s\circ h_2:\mathcal{SG}\to[0,1]$, given by
\[
f(x)=(s\circ h_2)(x)=\sin(h_2(x)),
\]
as in Example ~\ref{EX2}. Applying the framework in Corollary~\ref{main_theorem}, we can obtain new functions belonging to $\operatorname{dom}(\varepsilon)$. It is important to note, however, that the fact that $\operatorname{dom}(\varepsilon)$ is an algebra alone does not guarantee that all such newly constructed functions always belong to $\operatorname{dom}(\varepsilon)$.
\end{remark}

The following result provides a decomposition of finite-energy functions into nonlinear finite-energy components.

\begin{theorem}\label{sum_har}
Every nonnegative finite-energy function admits a nonlinear decomposition.
\[
f = f\left(\frac{1}{1+f} + \frac{f}{1+f} \right) = \frac{f}{1+f} + \frac{f^2}{1+f}.
\]
where both summands have finite energy.
\end{theorem}
\begin{proof}
By Corollary~\ref{main_theorem}, $\frac{f}{1+f}\in\mathrm{dom}(\varepsilon)$. Since $\mathrm{dom}(\varepsilon)$ is an algebra,
\[
\frac{f^2}{1+f} = f\cdot\frac{f}{1+f} \in \mathrm{dom}(\varepsilon).
\]
Hence, the required decomposition follows.
\end{proof}

\begin{theorem}\label{harm_decomp}
Any harmonic function that does not attain the value $-1$ on $\mathcal{SG}$ can be written as the sum of two finite-energy functions that are not harmonic.
\end{theorem}
\begin{proof}
Suppose $h$ is a non-constant harmonic function that does not attain the value $-1$ on $\mathcal{SG}$.
By Corollary ~\ref{coro_har}, it follows that $\frac{1}{1+h}$ and $\frac{h}{1+h}$ are in $\text{dom}(\varepsilon)$. Since $\text{dom}(\varepsilon)$ is an algebra, $h\cdot\frac{h}{1+h}=\frac{h^2}{1+h}\in\text{dom}(\varepsilon)$. Hence,
\begin{equation}\label{decomp_h}
h = h\left(\frac{1}{1+h} + \frac{h}{1+h} \right) = \frac{h}{1+h} + \frac{h^2}{1+h}.
\end{equation}\\
\textbf{Claim 1:} The function $\frac{h}{1+h}$ is not a harmonic function.\\
\textbf{Proof of claim 1:} Suppose $\frac{h}{1+h}$ is a harmonic function. For simplicity, let $h(P_1)=x, h(P_2)=y,$ and $h(P_3)=z$. Since $\frac{h}{1+h}$ is a harmonic function, it satisfies the harmonic extension algorithm. Hence,
\begin{align*}
\frac{h}{1+h}(P_{12})&=\frac{2}{5}\frac{h}{1+h}(P_{1})+\frac{2}{5}\frac{h}{1+h}(P_{2})+\frac{1}{5}\frac{h}{1+h}(P_{3})\\
\frac{h}{1+h}(P_{13})&=\frac{2}{5}\frac{h}{1+h}(P_{1})+\frac{2}{5}\frac{h}{1+h}(P_{3})+\frac{1}{5}\frac{h}{1+h}(P_{2})\\
\frac{h}{1+h}(P_{23})&=\frac{2}{5}\frac{h}{1+h}(P_{2})+\frac{2}{5}\frac{h}{1+h}(P_{3})+\frac{1}{5}\frac{h}{1+h}(P_{1}).
\end{align*}
Since $h$ is a harmonic function and $h(P_1)=x, h(P_2)=y,$ and $h(P_3)=z$, the above system of equations can be written as
\begin{align*}
\frac{\frac{2}{5}x+\frac{2}{5}y+\frac{1}{5}z}{1+\frac{2}{5}x+\frac{2}{5}y+\frac{1}{5}z}&=\frac{2}{5}\frac{x}{(1+x)}+\frac{2}{5}\frac{y}{(1+y)}+\frac{1}{5}\frac{z}{(1+z)}\\
\frac{\frac{2}{5}x+\frac{2}{5}z+\frac{1}{5}y}{1+\frac{2}{5}x+\frac{2}{5}z+\frac{1}{5}y}&=\frac{2}{5}\frac{x}{(1+x)}+\frac{2}{5}\frac{z}{(1+z)}+\frac{1}{5}\frac{y}{(1+y)}\\
\frac{\frac{2}{5}y+\frac{2}{5}z+\frac{1}{5}x}{1+\frac{2}{5}y+\frac{2}{5}z+\frac{1}{5}x}&=\frac{2}{5}\frac{y}{(1+y)}+\frac{2}{5}\frac{z}{(1+z)}+\frac{1}{5}\frac{x}{(1+x)}.\\
\end{align*}
Simplifying the above system of equations, it is obtained that
\begin{align*}
E_1:\quad x^2y+2x^2z+3x^2+xy^2-8xyz-4xy
+xz^2-2xz+2y^2z+3y^2+yz^2-2yz+2z^2=0,
\\[2mm]
E_2:\quad 2x^2y+x^2z+3x^2+xy^2-8xyz-2xy
+xz^2-4xz+y^2z+2y^2+2yz^2-2yz+3z^2=0,
\\[2mm]
E_3:\quad x^2y+x^2z+2x^2+2xy^2-8xyz-2xy
+2xz^2-2xz+y^2z+3y^2+yz^2-4yz+3z^2=0.
\end{align*}
Upon simplification of $E_1-E_2, E_2-E_3,$ and $E_3-E_1$, the following equations are obtained
\begin{equation}\label{Sys1}
\begin{aligned}
(y-z)(x^2+2x-yz-y-z) &= 0,\\
(x-y)(xy+x+y-z^2-2z) &= 0,\\
(x-z)(xz+x-y^2-2y+z) &= 0.
\end{aligned}
\end{equation}
CASE 1: Suppose $x\neq y\neq z$. Then, the system \eqref{Sys1} can be written as
\begin{align}
(x^2+2x-yz-y-z) = 0,\label{Sys2_eq1}\\
(xy+x+y-z^2-2z) = 0,\label{Sys2_eq2}\\
(xz+x-y^2-2y+z) = 0\label{Sys2_eq3}.
\end{align}
Eq.~\eqref{Sys2_eq3}$-$Eq. \eqref{Sys2_eq2} gives
\begin{equation*}
(z-y)(3+x+y+z)=0\implies 3+x+y+z=0 \text{ since } z\neq y.
\end{equation*}
Thus, 
\begin{equation}\label{sum_eqn}
y+z=-(x+3).
\end{equation}
Substituting Eq. \eqref{sum_eqn} in Eq. \eqref{Sys2_eq1}, 
\begin{equation}\label{pdt_eqn}
    yz=3+3x+x^2.
\end{equation}
By Eqs. \eqref{sum_eqn} and \eqref{pdt_eqn}, it can be concluded that $y$ and $z$ are the solutions of the quadratic equation
\begin{equation}\label{quad_eqn}
t^2+(x+3)t+(3+3x+x^2)=0.
\end{equation}
Since $y,z\in\mathbb{R}$, discriminant of Eq. \eqref{quad_eqn} should be non-negative. i.e., 
\begin{equation}\label{contra}
(x+3)^2-4(3+3x+x^2)\geq 0\implies -3(x+1)^2\geq 0
\end{equation}
 Eq.~\eqref{contra} leads to a contradiction, since the only possibility for Eq. \eqref{contra} to be true is $x=-1$, which is not true, since the harmonic function $h$ does not attain the value $-1$ at the boundary $P_1$.  \\
 CASE 2: Suppose $x=y\neq z$. Then, the system \eqref{Sys1} reduces to 
\begin{equation*}
(y-z)^2(1+y)=0\implies y=z \text{ or } y=-1,
\end{equation*}
both will, in turn, lead to a contradiction. The counterparts $x\neq y=z$ and $x=z\neq y$ follows similarly.
\par
It is obvious that $x=y=z$ is not possible, since $h$ under consideration is a non-constant harmonic function.\\
\textbf{Claim 2:} The function $\frac{h^2}{1+h}$ is not a harmonic function.\\
\textbf{Proof of Claim 2:} Suppose $\frac{h^2}{1+h}$ is a harmonic function. By Eq. \eqref{decomp_h}, it follows that \begin{equation*}
\frac{h}{1+h}=h-\frac{h^2}{1+h}
\end{equation*} is a harmonic function, since $h$ is a harmonic function. It is a contradiction to claim 1, and hence claim 2 follows.
\par
Hence, by the Eq.\eqref{decomp_h} along with claims 1 and 2, the result follows for non-constant harmonic functions.
\par
Now, suppose $h$ is a constant function, say $h(x)=c$. Let $h_0$ be any non-constant harmonic function that does not attain $-1$ on $\mathcal{SG}$.\\
\textbf{Claim 3:} The function $\frac{1}{1+h_0}$ is not a harmonic function.\\
\textbf{Proof of Claim 3:} Suppose $\frac{1}{1+h_0}$ is a harmonic function. Then, \begin{equation*}
\frac{h_0}{1+h_0}=1-\frac{1}{1+h_0}
\end{equation*} is a harmonic function. It is a contradiction to claim 1, and hence claim 3 follows.
\par
Now, \begin{equation*}
h(x)=c=\frac{c}{1+h_0(x)}+\frac{c\cdot h_0(x)}{1+h_0(x)},
\end{equation*}
with $\frac{c}{1+h_0(x)}$ and $\frac{c\cdot h_0(x)}{1+h_0(x)}$ being non-harmonic by claims 1 and 3, since $h_0$ is non-constant harmonic function.
\par
Hence, the theorem follows.
\end{proof}
\begin{remark}\label{Rem_for1over_h}
Suppose $h$ is a non-constant harmonic function which do not attain the value $0$ on $\mathcal{SG}$. Then,
\begin{equation*}
\frac{1}{h}=\frac{1}{1+(h-1)}
\end{equation*}
with $h-1$ is a non-constant harmonic function, which does not attain $-1$. Hence, by claim 3 from Theorem \ref{harm_decomp}, it follows that $\frac{1}{h}$ is not a harmonic function.
\end{remark}
Besides compositions with Lipschitz maps, finite-energy functions are also preserved under the symmetries of the Sierpiński gasket.

\begin{theorem}
Let $\mathcal{G}$ denote the group of all reflections and rotations on the Sierpiński gasket and let $g\in\mathcal{G}$. Then $f\circ g\in \mathrm{dom}(\varepsilon)$ if and only if $f\in \mathrm{dom}(\varepsilon)$.
\end{theorem}
\begin{proof}
The $m^{th}$ stage energy of the function $f\circ g$ is
\begin{align*}
\epsilon_m(f\circ g)=\Big(\frac{5}{3}\Big)^m\sum\limits_{x\sim_m y} \lvert f(g(x))-f(g(y))\rvert^2.
\end{align*}
Since $g$ is either a rotation or reflection, adjacency is preserved. Therefore,
\begin{align*}
\sum\limits_{x\sim_m y} \lvert f(g(x))-f(g(y))\rvert^2=\sum\limits_{x\sim_m y} \lvert f(x)-f(y)\rvert^2,
\end{align*}
which implies that $\varepsilon_m(f\circ g)= \varepsilon_m(f)$ for every $m\ge0$. Passing to the limit yields $\varepsilon(f\circ g)=\varepsilon(f)$. Thus, $\varepsilon(f)<\infty$ if and only if $\varepsilon(f\circ g)<\infty$, that is, $f\circ g\in \mathrm{dom}(\varepsilon)$ if and only if $f\in \mathrm{dom}(\varepsilon)$.
\end{proof} 
The finite energy functions, constructed with the help of Lipschitz maps as given
in Theorem \ref{first lemma}, need not be in $\text{dom}(\Delta)$. The details are provided in Proposition \ref{notindom} and
Corollary \ref{cor_notindom}.

\begin{proposition}\label{notindom}
Let $f\in\mathrm{dom}(\Delta)$ be a non-constant function with $f(x)\in\mathbb{R}^+\bigcup\{0\}$ for all $x\in\mathcal{SG}$. Then, $\frac{1}{1+f}$ and $\frac{f}{1+f}$ do not belong to $\mathrm{dom}(\Delta)$. 
\end{proposition}
\begin{proof}
Since $f\in\mathrm{dom}(\Delta)$, it follows that $1+f\in\mathrm{dom}(\Delta)$.  Suppose that $\frac1{1+f}\in\mathrm{dom}(\Delta).$ Then, $1=(1+f)\frac1{1+f}$  is the product of two non-constant functions in $\mathrm{dom}(\Delta)$, contradicting Proposition~\ref{whatnot}. Hence, $\frac1{1+f}\notin\mathrm{dom}(\Delta).$ Now, suppose that $\frac{f}{1+f}\in\mathrm{dom}(\Delta).$ Since $\frac1{1+f} = 1-\frac{f}{1+f}, $ it would follow that $\frac1{1+f}\in\mathrm{dom}(\Delta),$ contradicting the first part of the proof. Therefore, $\frac{f}{1+f}\notin\mathrm{dom}(\Delta).$ 
\end{proof}   
The following consequence is immediate.
\begin{cor}\label{cor_notindom}
Let $f\in\mathrm{dom}(\Delta)$, then $s\circ f$ need not be in $\mathrm{dom(\Delta)}$, where $s$ is a Lipschitz map.
\end{cor}

Although the domain of the Laplacian is not preserved under the transformations considered above, every function in $\mathrm{dom}(\Delta)$ can nevertheless be decomposed into functions lying outside $\mathrm{dom}(\Delta)$.

\begin{theorem}
Any non-constant $f\in \mathrm{dom}(\Delta)$ with $f(x)\in\mathbb{R}^+\bigcup\{0\}$ for all $x\in\mathcal{SG}$ can be written as the sum of two functions which are not in $\mathrm{dom}(\Delta)$.
\end{theorem}
\begin{proof}
Note that $f(x)$ can be written as $f(x)=\frac{f(x)}{1+f(x)}+\frac{f(x)^2}{1+f(x)}$ for any $x\in\mathcal{SG}$. By Proposition~\ref{notindom}, $\frac{f}{1+f}\not\in\mathrm{dom}(\Delta)$. Suppose that $\frac{f^2}{1+f}\in \mathrm{dom}(\Delta)$.  Then, $f-\frac{f^2}{1+f}=\frac{f}{1+f}\in\mathrm{dom}(\Delta)$, which contradicts Proposition~\ref{notindom}. Hence, $\frac{f^2}{1+f} \notin \mathrm{dom}(\Delta),$ and the proof is complete.
\end{proof}

\section{A condition for a FIF to have finite energy}\label{Section2}
Fractal interpolation functions provide a flexible framework for constructing continuous functions whose graphs exhibit fractal behaviour. Besides their geometric significance, understanding when such functions possess finite energy is important from the perspective of analysis on fractals, since finite-energy functions form the natural domain of the standard Dirichlet form. In this section, sufficient and necessary conditions are established for the fractal interpolation functions introduced in \cite{ri2020fractal} to belong to the energy domain. The construction of the fractal interpolation function (FIF) is briefly recalled below, as it forms the basis for the results established in this section. 

Let $h_i:\mathcal{SG}\subseteq \mathbb{R}^2\to \mathbb{R}$ be a non-constant harmonic functions and let $s_i:\mathbb{R}\to\mathbb{R}$ be continuous functions for $i=1,2,3$. Consider the iterated function system $\{w_i:i=1,2,3\}$ on $\mathcal{SG}\times \mathbb{R}\subseteq \mathbb{R}^2\times R$ is defined by
\begin{align}\label{EQ_intrpoln1}
w_i(x,y)=\begin{pmatrix}   u_i(x)\\F_i(x,y) \end{pmatrix}=\begin{pmatrix} \frac{1}{2}(x+P_i)\\s_i(y)+h_i(x) \end{pmatrix}.
\end{align}
Equivalently,
\begin{align}
    w_i(x^1,x^2,y)=\begin{pmatrix}
    u^1_i(x^1,x^2)\\u^2_i(x^1,x^2)\\F_i(x^1,x^2,y)
\end{pmatrix}=\begin{pmatrix}
    \frac{1}{2}(x^1+p^1_i)\\ \frac{1}{2}(x^2+p^2_i)\\s_i(y)+h_i(x^1,x^2)
\end{pmatrix},
\end{align}
where $x=(x^1,x^2)\in\mathcal{SG}$, $y\in\mathbb{R}$.
\par 
Let $\nu: V_1\subseteq \mathbb{R}^2\to \mathbb{R}$ denote a function on $V_1$, which can be seen as the interpolation values on the set $V_1$. Suppose,
\begin{equation}\label{EQ_intrpoln2}
    w_i(P_j,\nu(P_j))=\begin{pmatrix}
    \frac{1}{2}(P_i+P_j)\\ \nu(\frac{1}{2}(P_i+P_j))
\end{pmatrix} \text{ for } i,j\in\{1,2,3\}.
\end{equation}
Comparing Eq.s \eqref{EQ_intrpoln1} and \eqref{EQ_intrpoln2}, 
\begin{equation*}
    s_i(\nu(P_j))+h_i(P_j)=\nu(\frac{1}{2}(P_i+P_j)) \text{ for } i,j\in\{1,2,3\}.
\end{equation*}
In other words, the boundary values of the harmonic function $h_i$ are uniquely determined as
\begin{equation*}
    h_i(P_j)=\nu(\frac{1}{2}(P_i+P_j))-s_i(\nu(P_j)) \text{ for } i,j\in\{1,2,3\}.
\end{equation*}
\par
Consider the subspaces $C^*(\mathcal{SG})$ and $C^{**}(\mathcal{SG})$ of $C(\mathcal{SG})$ as follows:
\begin{align*}
& C^*(\mathcal{SG})=\{f\in C(\mathcal{SG}): f(p)=\nu(p), p\in V_0\}\\ 
&C^{**}(\mathcal{SG})=\{f\in C^*(\mathcal{SG}): f(p)=\nu(p), p\in V_1\}.
\end{align*}
Define $d_{C(\mathcal{SG})}(g_1,g_2):=\max\limits_{x\in \mathcal{SG}}\big\lvert g_1(x)-g_2(x)\big\rvert$ for all $g_1,g_2\in C(\mathcal{SG})$.
For $f\in C^*(\mathcal{SG})$ and $ x\in \mathcal{SG}$, define
$Tf(u_i(x)):=s_i(f(x))+h_i(x)$ for $i=1,2,3$. In other words, for  $x\in u_i(\mathcal{SG})$, $Tf(x)=s_i(f(u_i^{-1}(x)))+h_i(u_i^{-1}(x))$.

The following result, established in \cite{ri2020fractal}, shows that the operator $T$ maps the prescribed function space into itself.

\begin{lem}\cite{ri2020fractal}
   $ Tf\in C^{**}(\mathcal{SG})$ for all $f\in C^*(\mathcal{SG})$.
\end{lem}

The existence and uniqueness of the associated fractal interpolation function follow from the contraction property of the operator $T$.

\begin{theorem}\label{thm from pcf}\cite{ri2020fractal}
    Assume that $s_i$ are Matkowski contractions with the same function $\phi$. Then the operator $T$ is a Matkowski contraction (considered as $T:C^*(\mathcal{SG})\to C^*(\mathcal{SG})$). Hence, there is a unique continuous function $f:\mathcal{SG}\to \mathbb{R}$, which is a fixed point of $T$.
    In particular, $f|_{v_1}=\nu$ and $f(u_i(x))=s_i(f(x))+h_i(x)$ for all $i=1,2,3$, $x\in \mathcal{SG}$ and
    $G=w_1(G)\bigcup w_2(G)\bigcup w_3(G)$, where $G$ is the graph of $f$.
\end{theorem}
\begin{remark}
  Theorem \ref{Intrpoln_thm1}  is obtained as a special case of Theorem \ref{thm from pcf} for $n=1$ if $s_i:\mathbb{R}\to\mathbb{R}$ are taken to be the contraction maps defined by $s_i(y)=\alpha_i y$, with $0<|\alpha_i|<1$ for $i=1,2,3$.
\end{remark}
The following proposition provides the first step towards characterizing finite-energy fractal interpolation functions by relating the energy of the interpolant to that of the associated mappings $s_i\circ f$.

\begin{theorem}\label{first prop}
Assume that $s_i$ are Matkowski contractions with the same function $\phi$ and let $f$ be the fractal interpolation function as given in Theorem \ref{thm from pcf}. Then $f\in \mathrm{dom}(\varepsilon)$ iff $s_i\circ f\in \mathrm{dom}(\varepsilon)$ for $i=1,2,3$.
\end{theorem}
\begin{proof}
The $m^{th}$ stage energy of the fractal interpolation function $f$ is \\
\begin{align*}
\varepsilon_m(f)&=\Big(\frac{5}{3}\Big)^m\sum\limits_{x\sim_m y}\lvert f(x)-f(y)\rvert^2\\
&=\Big(\frac{5}{3}\Big)^m\sum\limits_{x\sim_m y}\lvert Tf(x)-Tf(y)\rvert^2\\
&= \Big(\frac{5}{3}\Big)^m\sum\limits_{\substack{x\sim_m y\\x,y\in u_i(SG)}}\sum\limits_{i=1}^3\Big[\big\lvert s_i\circ f(u_i^{-1}(x))+h_i(u_i^{-1}(x))\\&\hspace{3cm}-s_i\circ f(u_i^{-1}(y))-h_i(u_i^{-1}(y))\big\rvert^2\Big]\\
&\leq 2\Big(\frac{5}{3}\Big)^m\sum\limits_{\substack{x\sim_m y\\x,y\in u_i(SG)}}\sum\limits_{i=1}^3\Big[\big\lvert s_i\circ f(u_i^{-1}(x))-s_i\circ f(u_i^{-1}(y))\big\rvert^2\\&\hspace{3cm}+\big\lvert h_i(u_i^{-1}(x))-h_i(u_i^{-1}(y))\big\rvert^2\Big]\\
&=\frac{10}{3}\Big(\frac{5}{3}\Big)^{m-1}\sum\limits_{\substack{x\sim_{m-1} y\\x,y\in SG}}\sum\limits_{i=1}^3\big[\lvert s_i\circ f(x)-s_i\circ f(y)\rvert^2+\lvert h_i(x)-h_i(y)\rvert^2\big]\\
&=\frac{10}{3}\left[\sum\limits_{i=1}^3\varepsilon_{m-1}(s_i\circ f)+\sum\limits_{i=1}^3 \varepsilon_{m-1}(h_i)\right].
\end{align*}
As $m$ goes to infinity,  $\varepsilon(f)\leq \frac{10}{3}\sum\limits_{i=1}^3[\varepsilon(s_i\circ f)+\varepsilon(h_i)]$. Since harmonic functions have finite energy and $s_i\circ f\in \mathrm{dom}(\varepsilon) \implies f\in \mathrm{dom}(\varepsilon)$.
\par
Conversaly, if $f\in\text{dom}(\varepsilon)$, it easily follows that $s_i\circ f\in \text{dom}(\varepsilon)$, for all $i=1,2,3,$ since each $s_i$ is a Matkowski contraction and hence a Lipschitz map by Remark \ref{First_rem}.
\end{proof}



\begin{remark}
Replacing the harmonic functions $h_i$ by continuous functions $\phi_i$ and choosing $s_i(y)=d_i y$ for constants $d_i$ yields the linear fractal interpolation function considered in \cite{ruan2010fractal}. A more explicit criterion for finite energy in this special case was obtained in \cite{ruan2010fractal}. Theorem~\ref{first prop} extend this viewpoint by providing finite-energy criteria for the considerably more general class of nonlinear fractal interpolation functions introduced in \cite{ri2020fractal}.
\end{remark}

\section{Conclusion}

This paper develops a systematic framework for constructing finite-energy functions on the Sierpiński gasket beyond the classical classes of harmonic and piecewise harmonic functions. The construction is based on compositions with Lipschitz maps and various classes of contractions, and several properties of the resulting functions are established, including the min--max property, energy estimates, and their relationship with the domain of the Laplacian.

The framework is further applied to fractal interpolation functions on the Sierpiński gasket. A sufficient condition for finite energy is established, and, when the scaling functions are Banach contractions, a necessary and sufficient condition is obtained. These results extend the existing finite-energy constructions and provide a broader class of examples for further analysis on the Sierpiński gasket.

The present paper focuses primarily on the standard energy and the standard probability measure in the study of the Laplacian on the Sierpiński gasket. More generally, various energy forms along with their corresponding Laplacians and harmonic functions arise in the analysis on Sierpinski gasket. Investigating whether analogous results hold in these more general settings provides a natural and interesting direction for future research.

\backmatter

\bmhead{Acknowledgements}
The first author acknowledges the financial support provided from the Prime Minister’s Research Fellowship (PMRF ID: 3201486) during the course of this work. The author also gratefully acknowledges the Indian Institute of Technology Tirupati for providing an excellent research environment, where this work was carried out.

The authors sincerely thank the reviewers for their valuable comments and constructive suggestions, which have significantly contributed to improving the quality and clarity of the paper.

\section*{Declarations}
\begin{itemize}
\item Funding: The corresponding author was supported by Prime Minister's Research Fellowship.
\item Conflict of interest/Competing interests: The authors have no conflicts of interest to disclose.
\item Ethics approval and consent to participate: Not applicable.
\item Consent for publication: Not applicable.
\item Data availability: Not applicable.
\item Materials availability: Not applicable.
\item Code availability: Not applicable.
\item Author contributions: Harsha Gopalakrishnan – Conceptualisation, Methodology, Formal Analysis, and Original Manuscript Writing; Srijanani Anurag Prasad – Supervision, Validation, Review, and Editing of the Manuscript.
\end{itemize}

\bibliography{refer.bib}

@article{sahu2020box,
  title={On the box-counting dimension of graphs of harmonic functions on the {S}ierpi{\'n}ski gasket},
  author={Sahu, Abhilash and Priyadarshi, Amit},
  journal={Journal of Mathematical Analysis and Applications},
  volume={487},
  number={2},
  pages={124036},
  year={2020},
  publisher={Elsevier}
}

@article{ruan2010fractal,
  title={Fractal interpolation functions on post critically finite self-similar sets},
  author={Ruan, Huo-Jun},
  journal={Fractals},
  volume={18},
  number={01},
  pages={119--125},
  year={2010},
  publisher={World Scientific}
}

@article{ri2020fractal,
  title={Fractal functions on the {S}ierpinski Gasket},
  author={Ri, SongIl},
  journal={Chaos, Solitons \& Fractals},
  volume={138},
  pages={110142},
  year={2020},
  publisher={Elsevier}
}

@book{kigami2001analysis,
  author    = {Kigami, Jun},
  title     = {Analysis on Fractals},
  series    = {Cambridge Tracts in Mathematics},
  volume    = {143},
  publisher = {Cambridge University Press},
  address   = {Cambridge},
  year      = {2001}
}

@book{strichartz2006differential,
  author    = {Strichartz, Robert S.},
  title     = {Differential Equations on Fractals: A Tutorial},
  publisher = {Princeton University Press},
  address   = {Princeton, NJ},
  year      = {2006}
}

@article{rakotch1962note,
  title={A note on contractive mappings},
  author={Rakotch, Ephraim},
  journal={Proceedings of the American Mathematical Society},
  volume={13},
  number={3},
  pages={459--465},
  year={1962}
}

@article{rhoades1977comparison,
  title={A comparison of various definitions of contractive mappings},
  author={Rhoades, Billy E},
  journal={Transactions of the American Mathematical Society},
  volume={226},
  pages={257--290},
  year={1977}
}

@article{ri2011some,
  title={Some properties of fractal interpolation functions on {S}ierpinski gasket},
  author={Ri, Song-Gyong and Ruan, Huo-Jun},
  journal={Journal of Mathematical Analysis and Applications},
  volume={380},
  number={1},
  pages={313--322},
  year={2011},
  publisher={Elsevier}
}

@article{li2017energy,
  title={Energy and {L}aplacian of fractal interpolation functions},
  author={Li, Xiao-hui and Ruan, Huo-jun},
  journal={Applied Mathematics-A Journal of Chinese Universities},
  volume={32},
  pages={201--210},
  year={2017},
  publisher={Springer}
}

@article{kigami1989harmonic,
  title={A harmonic calculus on the {S}ierpinski spaces},
  author={Kigami, Jun},
  journal={Japan Journal of Applied Mathematics},
  volume={6},
  pages={259--290},
  year={1989},
  publisher={Springer}
}

@article{ben1999not,
  title={What is not in the domain of the {L}aplacian on {S}ierpinski gasket type fractals},
  author={Ben-Bassat, Oren and Strichartz, Robert S and Teplyaev, Alexander},
  journal={Journal of functional analysis},
  volume={166},
  number={2},
  pages={197--217},
  year={1999},
  publisher={Elsevier}
}

@article{jachymski2007nonlinear,
  title={Nonlinear contractive conditions: a comparison and related problems},
  author={Jachymski, Jacek and J{\'o}{\'z}wik, Izabela},
  journal={Banach center publications},
  volume={1},
  number={77},
  pages={123--146},
  year={2007}
}

@article{ccelik2008fractal,
  title={Fractal interpolation on the {S}ierpinski Gasket},
  author={{\c{C}}elik, Derya and Ko{\c{c}}ak, {\c{S}}ahin and {\"O}zdemir, Yunus},
  journal={Journal of mathematical analysis and applications},
  volume={337},
  number={1},
  pages={343--347},
  year={2008},
  publisher={Elsevier}
}

@article{de2013pcf,
  title={PCF self-similar sets and fractal interpolation},
  author={de Amo, Enrique and Carrillo, Manuel D{\'\i}az and S{\'a}nchez, Juan Fern{\'a}ndez},
  journal={Mathematics and Computers in Simulation},
  volume={92},
  pages={28--39},
  year={2013},
  publisher={Elsevier}
}

@article{kigami1993harmonic,
  title={Harmonic calculus on pcf self-similar sets},
  author={Kigami, Jun},
  journal={Transactions of the American Mathematical Society},
  volume={335},
  number={2},
  pages={721--755},
  year={1993}
}

@article{pasupathi2024very,
  title={A very general framework for fractal interpolation functions},
  author={Pasupathi, R and Miculescu, Radu},
  journal={Journal of Mathematical Analysis and Applications},
  volume={534},
  number={2},
  pages={128093},
  year={2024},
  publisher={Elsevier}
}

@article{barnsley1988fractals,
  title={Fractals everywhere academic press inc},
  author={Barnsley, M},
  journal={San Diego},
  volume={105},
  year={1988}
}

@inproceedings{sabot1997existence,
  title={Existence and uniqueness of diffusions on finitely ramified self-similar fractals},
  author={Sabot, Christophe},
  booktitle={Annales scientifiques de l'Ecole normale sup{\'e}rieure},
  volume={30},
  number={5},
  pages={605--673},
  year={1997}
}

\end{document}